\documentclass[conference]{IEEEtran}
\IEEEoverridecommandlockouts
\usepackage{cite}
\usepackage{amsmath,amssymb,amsfonts}
\usepackage{algorithm}
\usepackage{algorithmic}
\usepackage{graphicx}
\usepackage{textcomp}
\usepackage{xcolor}
\def\BibTeX{{\rm B\kern-.05em{\sc i\kern-.025em b}\kern-.08em
    T\kern-.1667em\lower.7ex\hbox{E}\kern-.125emX}}
\begin{document}

\title{Adaptive Memory Procedure for Solving Real-world Vehicle Routing Problem}
\author{
\IEEEauthorblockN{Nikica Perić\textsuperscript{a,*}, Slaven Begović\textsuperscript{b} and Vinko Lešić\textsuperscript{c}}
\IEEEauthorblockA{\textit{Laboratory for Renewable Energy Systems} \\
\textit{University of Zagreb, Faculty of Electrical Engineering and Computing} \\
Zagreb, Croatia, \\
\textsuperscript{a,*}nikica.peric@fer.hr,
\textsuperscript{b}slaven.begovic@fer.hr,
\textsuperscript{c}vinko.lesic@fer.hr
}}

\maketitle

\begin{abstract}
Logistics and transport are core of many industrial and business processes. One of the most promising segments in the field is optimisation of vehicle routes. Scientific effort is focused primarily on algorithms developed in simplified environment and cover a fraction of real industrial application due to complex combinatorial algorithms required to be promptly executed. In this paper, a real-world case study in all its complexity is observed and formulated as a real-world vehicle routing problem (VRP). To be able to computationally cope with the complexity, we propose a new procedure based on adaptive memory metaheuristic combined with local search. The initial solution is obtained with Clarke-Wright algorithm extended here by introducing a dropout factor to include a required stochastic attribute. The procedure and corresponding algorithms are tested on the existing benchmarks and further on the real industrial case study, which considers capacities, time windows, soft time windows, heterogeneous vehicles, dynamic fuel consumption, multi-trip delivery, crew skills, split delivery and, finally, time-dependent routes as the most significant factor. In comparison with the current state-of-the-art algorithms for vehicle routing problem with a large number of constraints, we obtain an average savings of 2.03\% in delivery time and 20.98\% in total delivery costs.
\end{abstract}


\section{Introduction}

Energy use for transport accounted for 30\% of global energy consumption in 2021, with road transport representing the bulk of the sector’s energy demand 78\% \cite{b46}. Also, the transport sector accounts for 10\% of EU GDP and any additional savings, even those less than 5\%, significantly contribute to the area. One way to achieve significant savings in this sector is the optimization of delivery logistics by Vehicle Routing Problem (VRP).

The VRP is a long-known combinatorial NP-hard problem where the goal is to allocate certain resources (packages) to a set of customers (delivery points) with the lowest cost. Computationally, it is among the most complex problems and for practical cases it is impossible to find global optimum in a limited time \cite{b48}. Therefore, various heuristic and metaheuristic methods \cite{b8, b18, b15, b4, b56} are developed to find a suboptimal or near-optimal solution.

The heuristics, also known as classical heuristics, is a collective name for strategies exploring promising areas of the search space. According to \cite{b6}, heuristics are divided into three main categories: constructive, two-phase, and improvement methods. Constructive heuristics are simple methods of low computational complexity that build a complete solution without an initial solution. 

Some of the well-known constructive heuristics are the sequential insertion algorithm \cite{b7} and the Clarke-Wright algorithm \cite{b8}, with, respectively, 10.42\% and 7.58\% of average suboptimality as reported in \cite{b44}. Both algorithms additonally inspired various improvements and many applications. 
Two-phase heuristics, in addition to principles similar to constructive heuristics, add different clustering strategies \cite{b9}. 
Improvement methods are based on local search, i.e. they are used to find a local optimum for the existing solution. Depending on the available computing time, these heuristics can be used to search for the optimum of a larger local neighborhood or to find the optimum in a smaller neighborhood in the shortest possible time. Improvement methods typically have a lower suboptimality, but longer execution time than constructive heuristics. The characteristic principle of improvement methods is the search of n-opt neighborhood introduced in \cite{b10}. In this way, an arrangement for $g_s$ selected vertices is found that minimizes the total cost. This is repeated for all combinations of vertices, giving $O^{g_s}$ complexity. Among the first significant applications of such heuristics, are \cite{b11} and \cite{b12}, where some of the local neighborhood search strategies are examined in detail.

An increase in computing power capabilities enabled computationally more intensive methods that provide less suboptimal solutions. The collective name is the metaheuristics, and they are the most commonly researched group of methods for solving a VRP \cite{b48}. The reason is the high accuracy, but also the large set of methods they include. Many metaheuristic algorithms use classical heuristic as base algorithms and then optimize their solutions in various ways. According to \cite{b13}, metaheuristics can be divided into single-solution-based and population-based. The most commonly used single-solution-based metaheuristics are tabu search \cite{b18} and variable neighborhood search \cite{b15}. Population-based metaheuristics are used less often, and the most popular is the genetic algorithm \cite{b16}.

In tabu search, the solution obtained by a local search (local optimum) is changed by finding a new, worse solution in the local neighborhood, that may take the solution to a new, potentially better local neighborhood. In order not to return to the same local optimum, the algorithm has a list of forbidden states, hence the name tabu. Well-known applications of tabu search are \cite{b17} and \cite{b18} where the authors found the best known solutions of 13/14 instances for the dataset from \cite{b44} in terms of the shortest path traveled. Among recent applications of the method is \cite{b53} where the authors solve two-echelon VRP with time windows and simultaneous pickup and delivery.

The concept of adaptive memory is introduced and combined with tabu search in \cite{b18}. Since this is a population-based algorithm that does not generate solutions by itself, it is combined with some single-solution-based algorithm, most often with tabu search. Adaptive memory is a pool of solutions with lower cost which dynamically updated during the search process. Adaptive memory in general characteristically introduces many degrees of freedom. This makes the algorithms more difficult to develop, but it also enables a more detailed search for solutions, which can make dramatic differences in the ability to solve problems \cite{b52}. Although this approach was used to find the first optimal solution to the problem from \cite{b44}, Adaptive Memory Procedure (AMP) rarely appears in the later literature.

In \cite{b19}, authors proposed the BoneRoute method, where sequences of vertices are referred to as bones. In \cite{b20}, new solutions of VRP are built by combining certain sequences of consecutive customers. Authors of \cite{b21} combine AMP and tabu search to solve fixed fleet open vehicle problem. To our knowledge, the hybrid of adaptive memory with the classical iterative memoryless method is solely used in \cite{b22}. Here the authors combine AMP and Variable Neighbourhood Search for solving real practical VRP.

The variable neighborhood search method is proposed in \cite{b15}. This method adds a so-called shaking step to the usual local search. Shaking step explores distant neighborhoods of the current existing solution and further uses the new solution regardless of its cost function. The whole procedure is repeated iteratively. One of the well-known applications of this algorithm solves the Periodic VRP \cite{b23}. In the Periodic VRP, the classical VRP is extended to a planning horizon of several days.

In addition to different methods of solving the VRP, there are also various real-application phenomena captured in the problem formulation, i.e. the variants. The variants respect practical apsects and further complicate the problem by introducing new constraints to the basic formulation of the VRP. Despite the fact that real-world problems often require a large number of variants of VRP, most of the literature cover only the basic ones. According to \cite{b13}, Capacitated VRP is applied in 98.91\% of relevant journal papers with metaheuristics, VRP with Time Windows in 37.32\%, VRP with Pickup and Delivery in 16.3\%, and Heterogeneous VRP in 13.41\% of papers. Due to introduction of additional complexity, other variants are used in less than 10\% of the literature. This opens space for the development of new metaheuristic methods to respect more variants, individually or jointly.

The variants considered in this paper are:

\begin{itemize}
\item[--] Capacitated VRP (CVRP),
\item[--] VRP with Time Windows (VRPTW),
\item[--] VRP with Soft Time Windows (VRPSTW),
\item[--] Heterogeneous VRP (HVRP),
\item[--] Multi-Trip VRP (MTVRP),
\item[--] Skill VRP (SVRP),
\item[--] Split Delivery VRP (SDVRP),
\item[--] Time-Dependent VRP (TDVRP),
\end{itemize}

\noindent{with various additional costs and cost functions.}

Most established variants are CVRP, introduced in \cite{b24}, and VRPTW, elaborated in \cite{b25}. The CVRP considers the mass or volume restrictions of the vehicle, jointly referred to as the capacity, and the amount of demand ordered, which creates the need to return the vehicle to the warehouse. This is different from traveling salesman problem, where no location needs to be reached more than once. The VRPTW introduces time windows as specific time slots in which demands can be delivered. The variant respects staff working hours. Authors in \cite{b43} found 42 out of 56 optimal solutions in dataset from \cite{b25}, the most common VRPTW benchmark. 
The HVRP \cite{b26}, \cite{b22} enables the use of vehicles of different capacities within the same problem. In this way, the fleet of vehicles can be distributed in several different ways with regard to the geographical locations of the delivery points and their time windows, which results in a reduction of the total cost of delivery.
The MTVRP \cite{b27}, \cite{b41} is applied in addition to the VRPTW in order to fill the vehicle time windows as much as possible and thus reduce the number of vehicles required. This variant is especially important for last-mile delivery as vehicles often return to the depot several times during the working day. 
A variant that takes into account different vehicle characteristics and driver capabilities is SVRP \cite{b29}. This possibility is important in problems that consider both city and long-distance deliveries. For example, oversized vehicles cannot drive through narrow city streets, and smaller vehicles are necessary. Also, drivers may have different roles within the company, and some of them cannot drive to any delivery point.
The SDVRP \cite{b30}, \cite{b47} enables order splitting into several smaller ones, i.e. an order can be delivered by visiting a location more than once. In this way, better utilization of vehicle space can be achieved. Also, by adding this constraint and combining it with the bin packing problem \cite{b31}, it is possible to significantly improve solutions for cases with packages of various shapes.
Demands and arrival times can vary considerably due to different day of the week, increased traffic congestion or some other specific situation. Therefore, it is beneficial to allow vehicles to arrive outside the ideal time window with a certain penalty. Similarly, overtime for vehicles is allowed. This is respected by VRPSTW, introduced in \cite{b32}.

In the paper, the largest focus is put on TDVRP. Although, according to \cite{b13}, this variant is considered in only 3.62\% of papers due to high complexity, the authors consider it as one of the key variants for solving real-world problems. In our research, VRP without time-dependent constraints produced between 3\% and 154\% travel time mismatch. The TDVRP is first mentioned in \cite{b33}. 
In \cite{b34}, authors recognize that optimal solutions without TDVRP are mostly suboptimal in real applications, and are often infeasible for TDVRP. On the dataset they proposed, for less variable driving durations, TDVRP achieved 0.9\% time savings, and 1.18\% of the solution without TDVRP was infeasible. For more variable driving durations, the time saving is 15.4\%, and the share of infeasible solutions without TDVRP is 39.3\%. Authors define TDVRP as a variant in which a new dimension is added to the standard two-dimensional time matrix, depending on the division of the time horizon into intervals. 
In \cite{b35}, the combination of TDVRP and VRPTW is solved using the Ant Colony System algorithm hybridized with insertion heuristics. It is emphasized that the TDVRP variant is particularly important in an urban context where the traffic plays a significant role. Most articles that analyze TDVRP divide traffic congestion into two categories. First is predictable, such as the heavy traffic during the daily peak hours, and other is less predictable, caused by accidents and weather conditions. The authors in \cite{b36} state that 70-87\% of all traffic congestion delays belong to the predictable category, which are also easier to solve, so most papers are based only on them. The authors also consider the non-passing property, which means that the vehicle must reach the destination earlier if it started driving earlier. The authors in \cite{b54} have defined time-dependent speed functions at the network level to deal with it. In \cite{b37} and \cite{b55}, the TDVRP is combined with the Green VRP, where electric vehicles are used. The problem is solved using mixed integer linear programming, followed by a combination of clustering, routing and simulated annealing. Authors in \cite{b42} combine sweep algorithm and improved particle swarm optimization to solve TDVRP and VRPSTW with stohastic factors.

The paper expands the current state of the art with three contributions. First contribution is the addition of a dropout factor to the widely used Clarke-Wright algorithm \cite{b8}. This adds stochasticity to the algorithm, which improves its application as an initial solution in metaheuristics. We refer to it as Randomized Constrained Clarke-Wright algorithm.

Second contribution is the combination of local search method for intensification and Adaptive memory procedure for diversification of the search process.

Many variants of VRP have been described, but the works in literature usually solve one or a few of them, which is not satisfactory for real-world problems. Third contribution of this paper is therefore the derived VRP as a combination of the variants listed in the introduction: CVRP, VRPTW, VRPSTW, HVRP, MTVRP, SVRP, SDVRP, TDVRP. 
We could not find a combination of these variants in the literature. Moreover, some of them are rarely listed even independently. Together, first two contributions enable the fast AMP to promptly and efficiently execute highly complex problems of many VRP variants. They are validated on standard benchmark and real world case study.

The paper is organized as follows. The mathematical formulation of the problem with all constraints is described in Section II. Section III provides a detailed description of the algorithm. In Section IV, existing benchmark problem and our real-world case are described, and Section V shows results in comparison to the most common state of the art algorithms. At the end, the conclusion is given in Section VI.

\section{Mathematical Formulation}

The VRP can be described as a problem of solving the following graph. Let $G_0 = (V_0, A_0)$ be a complete graph, where $V_0 = \{0, ..., I^0\}$ denotes the vertex set and $A_0 = \{(i,j):i,j\in V_0,i\neq j\}$ is the arc set with notation explained in Table~I. A vertex in VRP is often referred to as the location. Vertex $i=0$ corresponds to the depot (starting point), and others are related to customers (delivery points).

The vehicle fleet is heterogeneous, which means it contains $N$ vehicles of different capacities $Q_n$. Every delivery point has a capacity demand $q_i$ and a condition to fit the demand into a single vehicles is:

\begin{equation}
0 \leq q_i \leq Q^{max}.
\end{equation}

In real-world problems, some demands are often larger than the maximum capacity of a single vehicle (SDVRP). In this case, the demands, and thus the vertices of the graph, can be separated into smaller parts that satisfy (1) and the graph is expanded with new vertices and arcs. The newly created graph consists of vertex set $V = \{0, ..., I\}$, where $I \geq I^0$, with corresponding arc set $A$. Vertex $i=0$ is still the depot, and others are customers.

In TDVRP, the travel time between two locations is denoted by $t_{ijk}$. Given that the travel durations are predictions for one or more time sections $k \in K$, they satisfy non-passing property. The travel time between $i$ and $j$ is given by prediction for time section $k$ regardless of the duration of the ride.

The main parameter of the cost function is the travel cost. In simple options, it is usually assumed that travel cost is equal to the distance between two locations ($c_{ij}=d_{ij}$) or travel time between two locations ($c_{ijk}=t_{ijk}$ for TDVRP).

In more complex options, monetary cost is often calculated. In this paper, the monetary cost is calculated as the product of the price of fuel and the amount of fuel consumed. Fuel consumption is calculated using the equation:
\begin{equation}
c_{nijk} = (c_{n}^0 + \mu_n^m \cdot m_{ijk}) \cdot (1 - \mu_n^v \cdot v_{ijk}),
\end{equation}
where constants $c_{n}^0$, $k_n^m$ i $k_n^v$ are characteristic for the type of vehicle and driver, $m_{ijk}$ is the mass of the current load, and $v_{ijk}$ is the average speed on the way from location $i$ to $j$ in time section $k$.

Binary decision variable $x_{nijk}$ gives information if the vehicle $n$ traverses an arc $(i, j)$ within a time section $k$ in found solution ($x_{nijk}=1$) or not ($x_{nijk}=0$). 
Constraints on the variable $x_{nijk}$ apply if vehicle $n$ cannot visit location $i$ (SVRP). Then we can assume:

\begin{equation}
x_{nijk} = 0,    \quad \quad \quad \quad \forall j, k,
\end{equation}
\begin{equation}
x_{njik} = 0,    \quad \quad \quad \quad \forall j, k.
\end{equation}

In addition, there are constraints depending on the time window in which it is acceptable to arrive at locations or use a vehicle (VRPTW). Each location $i$ implies the interval when it can be visited:

\begin{equation}
[b_i^{vr}, e_i^{vr}], \quad b_i^{vr} < e_i^{vr}.
\end{equation}

Every vehicle $n$ is allowed only in the interval

\begin{equation}
[b_n^{vh}, e_n^{vh}], \quad b_n^{vh} < e_n^{vh}.
\end{equation}

In particular, $x_{nijk}=0$ if whole time section $k$ is out of interval $[b_i^{vr}, e_i^{vr}]$ for location $i$. The same applies to location $j$. Also, $x_{nijk}=0$ if time section $k$ is out of interval $[b_n^{vr}, e_n^{vr}]$ for vehicle $n$. Soft time windows, denoted as uppercase index $s$, are also introduced for locations and vehicles:

\begin{equation}
[b_i^{vrs}, e_i^{vrs}], \quad b_i^{vr} \leq b_i^{vrs} < e_i^{vrs} \leq e_i^{vr},
\end{equation}

\begin{equation}
[b_n^{vhs}, e_n^{vhs}], \quad b_n^{vh} \leq b_n^{vhs} < e_n^{vhs} \leq e_n^{vh}.
\end{equation}

Using soft time windows imply adding penalty of exceeding the designated arrival time. The penalty is added outside intervals from (7) and (8). Three parameters are defined for each early arrival at the location: $a_0^b$, $a_1^b$ multiplied by the pre-arrival time, and $a_2^b$ multiplied by the time the vehicle spent at the location outside the desired time window. 
Analogously, for late departures from the location, the parameters $a_0^e$, $a_1^e$ and $a_2^e$ apply. If the arrival time of the vehicle to location as $b_i$, and departure time $e_i$ is defined, then the penalty cost of early arrival $w$ to location $i$ is:

\begin{equation}
w_{i}^{ve,b} = a_0^b + a_1^b \cdot (b_i^{vrs} - b_i) + a_2^b \cdot (\min(e_i, b_i^{vrs}) - b_i),
\end{equation}

\noindent{and the penalty cost for late arrival at the location is:}

\begin{equation}
w_{i}^{ve,e} = a_0^e + a_1^e \cdot (e_i - e_i^{vrs}) + a_2^e \cdot (e_i - \max(b_i, e_i^{vrs})).
\end{equation}

The analogy is applied for the vehicle time windows calculation of corresponding $w_{n}^{vh,b}$ and $w_{n}^{vh,e}$.

The goal of the developed VRP algorithm is to minimize:


\begin{equation}
\label{eq:1}
\begin{aligned}
J=&\sum_{n=1}^{N}\sum_{i=0}^{I}\sum_{j=0}^{I}\sum_{k=1}^{K} c_{nijk} \cdot x_{nijk} + \\
&+\sum_{i=0}^{I} (w_{i}^{ve,b} + w_{i}^{ve,e}) + 
\sum_{n=1}^{N} (w_{n}^{vh,b} + w_{n}^{vh,e}),
\end{aligned}
\end{equation}

\noindent{s.t.}

\begin{equation}
\sum_{n=1}^{N}\sum_{j=0}^{I}\sum_{k=1}^{K}x_{nijk} = 1 \quad \quad \quad \quad \quad \quad \forall i>0,
\end{equation}
\begin{equation}
Q_{n^*} \geq \sum_{i=1}^{I}\sum_{j=0}^{I}\sum_{k=1}^{K}q_i x_{n^*ijk} \quad \quad \forall n^*.
\end{equation}

The variables used in the formulation are listed in Table~I. In this paper, travel time and monetary cost are both taken as the travel cost, which is described in more detail in Section~V.

In (13), $n^* \in \{1,...,N^*\}$ is introduced instead of $n$ to enable MTVRP, and represents one drive. In all drives, all delivery points must be visited once, as well as the vehicles in (12). The locations visited in one drive are a subset of the locations visited by one of the vehicles. The drive consists of all delivery points visited between the two arrivals of the vehicle at the depot. Therefore, each used vehicle performs 1 or more drives.

\begin{table}[htbp]
\caption{Variables used in the formulation}
\begin{center}
\begin{tabular}{|c|c|}
\hline
\multicolumn{2}{|c|}{\textbf{Indices}} \\
\hline
$i$ & Index of the vertex (location), arc start; $i \in \{0,...,I\}$        \\
$i=0$ & Depot (starting point)      \\
$i \in \{1,...,I\}$ & Customer (delivery point)       \\
$j$ & Index of the vertex (location), arc end; $j \in \{0,...,I\}$    \\
$k$ & Time section number; $k \in \{1,...,K\}$    \\
$n$ & Vehicle number; $n \in \{1,...,N\}$        \\
$n^*$ & Number of drive; $n^* \in \{1,...,N^*\}$        \\
\hline
\multicolumn{2}{|c|}{\textbf{Parameters}} \\
\hline
$I^0$ & Number of delivery points before split delivery        \\
$I$ & Number of delivery points after split delivery        \\
$Q_n$ & Capacity of vehicle $n$      \\
$Q^{max}$ & Capacity of the vehicle with largest capacity      \\
$q_i$ & Demand of delivery point $i$ (capacity)     \\
$d_{ij}$ & Road distance from location $i$ to location $j$     \\
$t_{ijk}$ & Travel time from location $i$ to location $j$ in time section $k$     \\
$b_i$, $e_i$ & Arrival time to and departure time from location $i$  \\
$b_i^{vr}$, $e_i^{vr}$ & Time window start and end for location $i$  \\
$b_n^{vh}$, $e_n^{vh}$ & Time window start and end for vehicle $n$  \\
$b_i^{vrs}$, $e_i^{vrs}$ & Soft time window start and end for location $i$  \\
$b_n^{vhs}$, $e_n^{vhs}$ & Soft time window start and end for vehicle $n$  \\
$x_{nijk}=1$ & Vehicle $n$ traverses an arc $(i,j)$ in time section $k$      \\
\hline
\multicolumn{2}{|c|}{\textbf{Costs and penalties}} \\
\hline
$c_{nijk}$ & Travel cost of vehicle $n$ from location $i$ to location $j$ \\
&  in time section $k$     \\
$m_{ijk}$ & the mass of the current load \\
$v_{ijk}$ & average speed from location $i$ to $j$ in time section $k$ \\
$a_0^{b}$, $a_0^{e}$ & Early and late delivery penalty cost  \\
$a_1^{b}$, $a_1^{e}$ & Early and late delivery penalty cost failed time factor  \\
$a_2^{b}$, $a_2^{e}$ & Early and late delivery penalty cost duration factor  \\
$w_i^{ve,b}$, $w_i^{ve,e}$ & Early and late delivery penalty cost for delivery point $i$  \\
$w_n^{vh,b}$, $w_n^{vh,e}$ & Early and late delivery penalty cost for vehicle $n$  \\
\hline
\end{tabular}
\label{tab1}
\end{center}
\end{table}

\section{Algorithm design}

Here the AMP is used to solve the VRP. It consists of three parts:

\begin{itemize}
\item[--] Initial solution generation,
\item[--] Local search,
\item[--] Adaptive memory procedure. 
\end{itemize}

\subsection{Initial solution generation}

\begin{figure*}[!h]
\centering
\includegraphics[width=.98\textwidth]{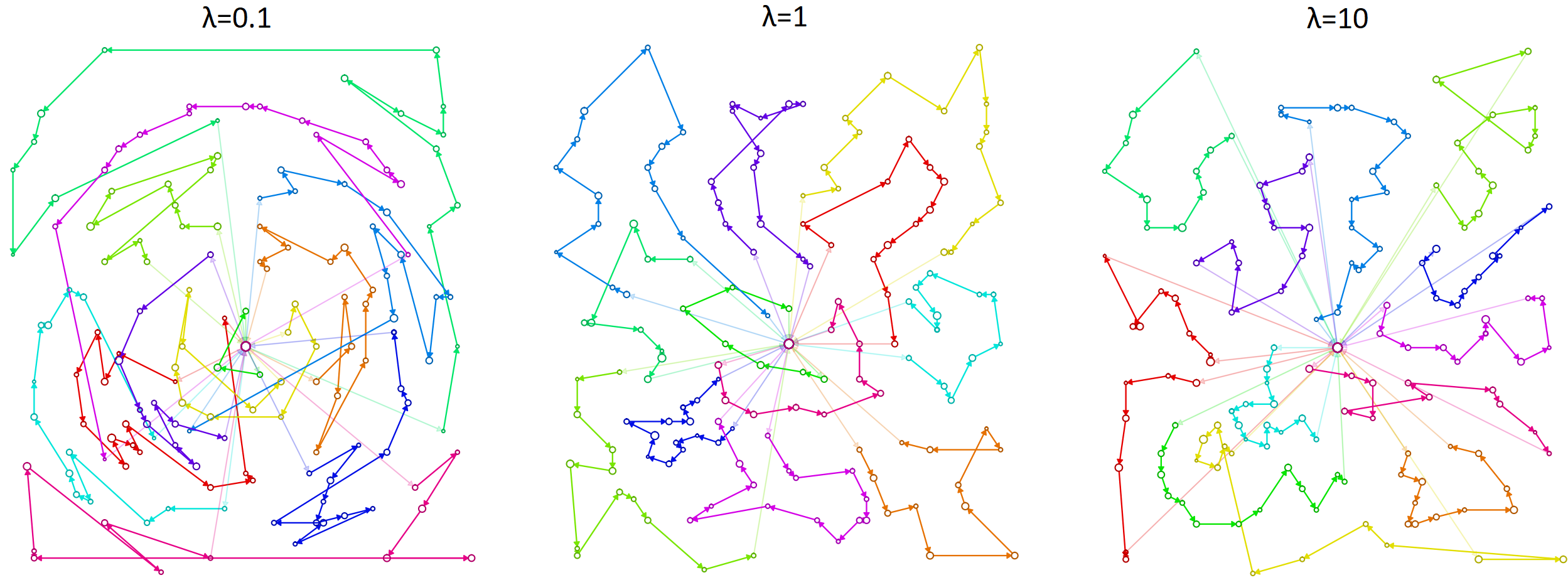}
\caption{Comparison of CVRP solutions using extended Clarke-Wright algorithm with different $\lambda$ values}
\label{FIG:hc}
\end{figure*}

As local search generates solution by improving the initial solution, it takes longer to improve the solution of a higher initial cost. The premise is that for the initial solution of a sufficiently low cost the algorithm converges faster and provides the final solution of a lower suboptimality. Therefore, large attention is given to the spatial layout and other parameters during route generation. For this reason, the initial solution is most often generated by constructive heuristics. One of them is the Clarke-Wright algorithm \cite{b8} with a wide adoption credited to fast execution. The original version of the algorithm first builds a solution where each of the $I$ delivery points has its own route (0, $i$, 0), and then these routes are interconnected through an iterative process. The connection order is determined by savings matrix $S$, where matrix elements are calculated as:

\begin{equation}
s_{ij} = c_{i0} + c_{0j} - c_{ij}.
\end{equation}

Routes with a higher $s_{ij}$ are joined first followed by routes with a lower one. Routes are connected in such a way that (0, ..., $i$, 0) and (0, $j$, ..., 0) merge into the route (0, ..., $i$, $j$, ..., 0). This is repeated iteratively as long as the capacity constraints allow it. The algorithm is deterministic, which achieves a high ratio of accuracy and low computational complexity, but also limits the set of solutions that can be obtained in case this algorithm is used only as an initial solution.

Since AMP requires several different initial solutions, it is beneficial to add stochasticity to the initial solution. We achieve it in two ways. The first way is to add the route shape parameter $\lambda$ to the savings calculation:

\begin{equation}
s_{ij} = c_{i0} + c_{0j} - \lambda \cdot c_{ij}.
\end{equation}

This parameter is proposed in \cite{b38} because the original Clarke and Wright algorithm tends to produce good routes at the beginning of execution, but worse routes toward the end. Depending on this parameter, the routes have a straight or rounded shape, as shown in Fig.~1.
Values of $\lambda$ in the figure are chosen to illustrate the influence on solution shape. Values that give good results are usually between 0.4 and 2 \cite{b45}.

Although $\lambda$ makes it possible to find several different solutions, the number of solutions is still relatively small and there are parts of solutions that are repeated regardless of the value of this parameter. Such pieces are usually part of the optimal solution but sometimes they prevent the algorithm from finding the optimal solution because it is difficult to divide them through the further process of local search. To avoid such pieces of solutions constantly directing the optimization to the same local optimum, we propose to introduce a dropout parameter, denoted as $p_d$. This parameter determines a percentage of the savings matrix ignored when generating the initial solution. A higher $p_d$ increases the average cost but also increases the diversity of the obtained initial solutions. 

This effect is analyzed in more detail in Fig.~2, where the influence of different values of parameters $\lambda$ and $p_d$ is shown.
For each combination of these parameters, the generation of the initial solution was run 15 times on 5 different instances from \cite{b44}, which is a widely adopted established benchmark introducing Christofides, Mingozzi, Toth (CMT) instances. Solution diversity represents the ratio of different arcs obtained in 15 solutions and the number of arcs required for one solution. Average cost is calculated over all 5 instances and 15 runs. Changing $\lambda$ in the range 0.2-1.8 increases solution diversity while maintaining a slightly lower average cost. Introducing $p_d$ has the same trade-off, but can more aggressively affect the initial solution. By combining these two parameters, the initial solution can be adapted to the given requirements. Thus, for search diversification, it is beneficial to increase $p_d$ and the search range of $\lambda$, while for search intensification, it is beneficial to reduce these two parameters. 
The influence of these two parameters was examined in the next chapter after starting the local search.

\begin{figure}[!h]
\centering
\includegraphics[width=.48\textwidth]{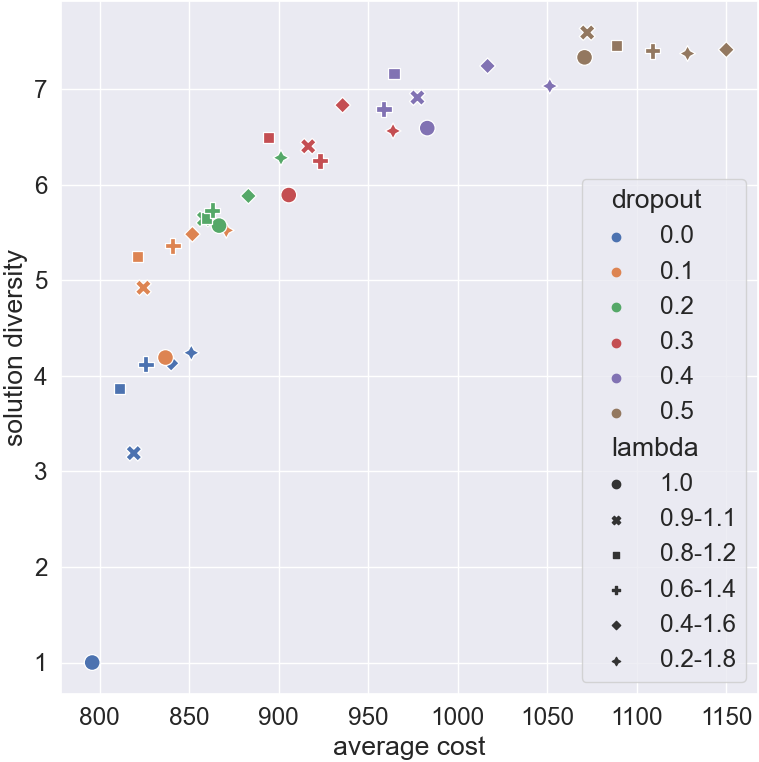}
\caption{Average cost and solution diversity of the initial solution depending on the route shape parameter $\lambda$ and the dropout parameter $p_d$}
\label{FIG:hc}
\end{figure}

The effect of $p_d$ parameter is also shown on a specific example in Fig.~3. The figure compares two initial solutions, on the left without dropout, and on the right with dropout, both with $\lambda=1$. Cost with dropout is higher as expected, but a different prominent part of the solution is obtained. This allows the local search to search in the environment of a better local optimum, which finally results in increased accuracy.

\begin{figure*}[!h]
\centering
\includegraphics[width=.98\textwidth]{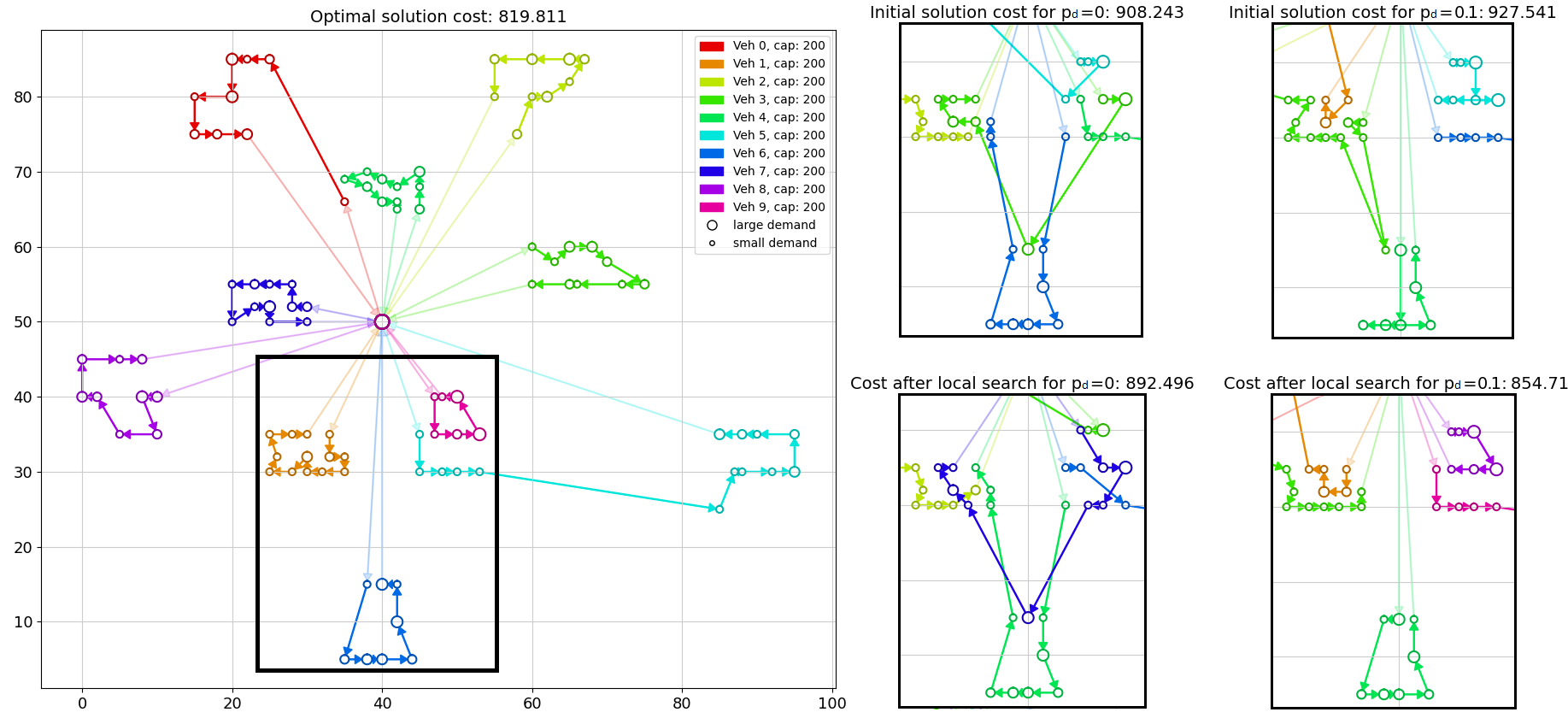}
\caption{Section of the CVRP solution after using the RCCW algorithm and local search without dropout and with dropout}
\label{FIG:hc}
\end{figure*}

The original Clarke-Wright algorithm is created for basic CVRP, which means that two routes with the largest $s_{ij}$ will not merge only if the sum of their capacities exceeds the vehicle capacity. Also, in the original algorithm it is not necessary to consider $s_{ji}$. Here, we improve the algorithm by introducing asymmetry of $s_{ij}$ and $s_{ji}$, adding them separately in savings matrix $S$.

To generate the initial solution, the Clarke-Wright algorithm is used with the addition of parameters $\lambda$ and $p_d$ and additional constraints. We refer to the algorithm as Randomized Constrained Clarke-Wright (RCCW). The pseudo-code of the algorithm is shown in Algorithm \ref{alg1}. Parameters $\lambda$ and $p_d$ are therefore used for tuning the algorithm. It turned out that in the initial phase of AMP, the better values are $0.1 \leq \lambda \leq 1.7$ and $p_d \in \{0.3, 0.8\}$ which bring more stochasticity. In this way, a more diverse set of solutions is found with a not too high cost in the initial phase of the algorithm, as shown in Fig.~2. In the later phase, it is beneficial to intensify the search with values $0.6 \leq \lambda \leq 1.2$ and $p_d \in \{0.1, 0.15\}$.

\begin{algorithm}
\caption{Randomized Constrained Clarke-Wright}
\label{alg1}
\begin{algorithmic}[1]
\STATE load inputs: vehicles, locations
\STATE select route shape $\lambda$
\STATE select dropout percentage $p_d$
\STATE generate solution of $I$ vehicle routes (0, $i$, 0)
\STATE calculate savings matrix $S$ of items $s_{ij}$ using (11)
\STATE sort values in savings matrix $S$ in decreasing order
\STATE drop $p_d$ items from $S$

\FOR {$s_{ij}$ in $S$}
    \IF {routes (0, ..., $i$, 0) and (0, $j$, ..., 0) can be merged}
        \STATE merge (0, ..., $i$, 0) and (0, $j$, ..., 0) to (0, ..., $i$, $j$, ..., 0)
    \ENDIF
\ENDFOR
\STATE compress solution to $N$ vehicles
\RETURN generated initial solution

\end{algorithmic}
\end{algorithm}

\subsection{Local search}

Local search is a heuristic for finding the local optimum of the existing solution (intensification of the search process). By increasing the number of strategies and the number of dimensions in which they are applied, the probability of an optimal solution showing up is significantly higher \cite{b10}, but it also slows down the search. The developed algorithm strives more for accuracy than speed of execution. The following search strategies were used in the algorithm: swap vertices, move vertices, swap drives, move drives and change drive direction.

Swap vertices works as follows. The vertices are divided into groups of size $g_s \in \{1, 2, ...\}$ and then the arrangement of the groups is shuffled to find the one that minimizes the cost. The number of different arrangements is $g_n!$, where $g_n \in \{2, 3, ...\}$ is the number of groups. An example for $g_n=2$ and $g_s=1$ is shown in Fig.~4 on the left.
In the move vertices strategy, only a single group is selected with one or more vertices. Additionally, the location in the solution where the group will be moved is selected. An example of this strategy for $g_s=1$ is shown in Fig.~4 on the right.
By using swap vertices and move vertices, vertices that were in the inappropriate drive in the initial solution are moved. The appropriate drive is the one where delivery to these vertices will add less impact to the total cost. In one iteration, many combinations of groups are usually tested and the one that reduces the value of the cost function the most is selected. By increasing the number of tested combinations, the complexity also increases exponentionally.
Therefore, in the first part of the algorithm, the variable $g_c$ is introduced, which determines how many different combinations with the selected $g_s$ and $g_n$ are tested.

\begin{figure*}[!h]
\centering
\includegraphics[width=.98\textwidth]{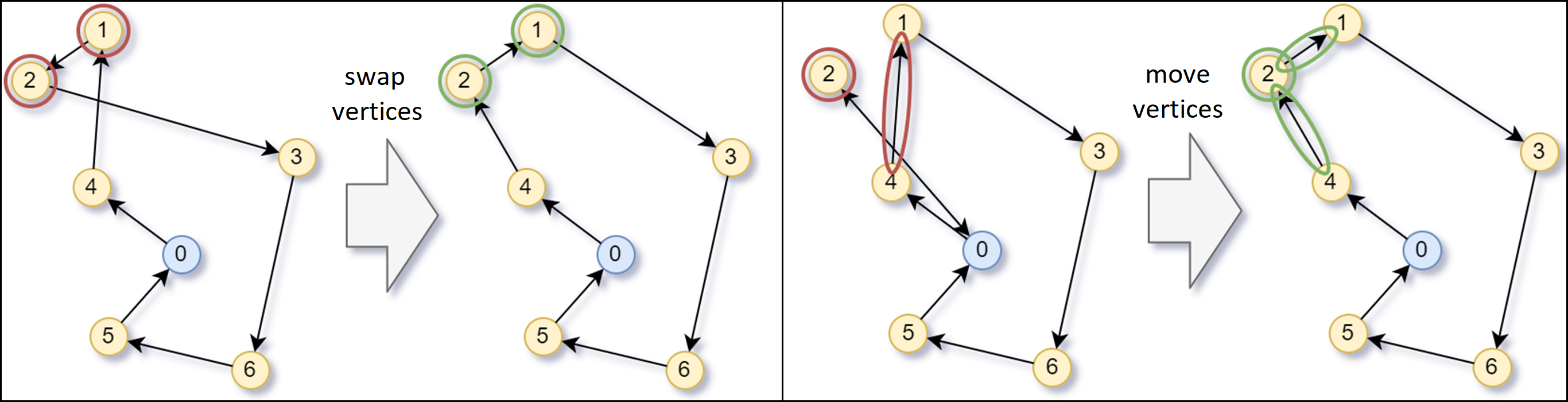}
\caption{Swap vertices and move vertices example}
\label{FIG:hc}
\end{figure*}

Swap drives and move drives use the entire drive instead of smaller groups. In this way, the drives are better adapted to time parameters (beneficial for VRPTW, VRPSTW and TDVRP) and vehicle parameters (beneficial for HVRP).
Change drive direction is applied on each drive to determine the direction in which the drive adds less value to the cost function (beneficial for TDVRP, or more generally for asymmetric VRP).

The above strategies are performed iteratively as long as the value of the cost function decreases. The pseudo-code of the local search algorithm is shown in Algorithm \ref{alg2}.

The strategy for each iteration is selected by the parameter $s$ which depends on the iteration number. The goal is that in the first part of the algorithm, where $g_c$ is relatively small, faster and simpler strategies are performed more often. Slower strategies, which search a larger local neighborhood, appear less often in the first part. Thus, swap nodes with $g_n$ == 3 and $g_s$ == 1 appear in 3\% iterations, swap nodes with $g_n$ == 2 and $g_s$ == 2 in 0.3\% iterations, swap nodes with $g_n$ == 2 and $g_s$ == 3 in 0.5\% iterations, and swap nodes with $g_n$ == 2 and $g_s$ == 1 in 47.5\% iterations.
Move node appears in 47.8\% of cases, while swap drives, move drive and find best are used in 0.3\% of iterations each.
In the second part of the local search, the maximum is selected for $g_c$, i.e. the entire local neighborhood is searched. In this phase, each strategy is repeated once as long as any of them can improve the current solution, i.e. until the local optimum is found.
Once a local optimum is found, these strategies are no longer beneficial, and further progress of the solution is achieved using AMP.

The result of the optimization after the first part of the local search is analyzed with different values of the parameters $\lambda$ and $p_d$. As for the initial solutions, the influence of $\lambda$ and $p_d$ was analyzed on 15 runs of each of the 5 instances. The test results are on Fig.~5. The advantages of a higher dropout and a larger range of $\lambda$ values are highlighted here. These values resulted in a worse average cost, but also a greater solution diversity in the initial solution. In particular, the case denoted as green circle in Fig. 2 (pd =, lambda =) resulted in the same green circle notation in Fig. 5, compared to a blue dot case of initial lowest cost but final highest cost. After the local search, such parameters brought a lower cost because the greater solution diversity increased the search space. Thus, the 2 lowest average costs were achieved with a dropout value of 0.2, and 3 of the 5 lowest costs with $0.2 \leq p_d \leq 0.4$ and $0.4 \leq \lambda \leq 1.6$. Based on this research, random values from the following intervals were selected: $0.2 \leq p_d \leq 0.4$ and $0.4 \leq \lambda \leq 1.6$.

\begin{figure}[!h]
\centering
\includegraphics[width=.48\textwidth]{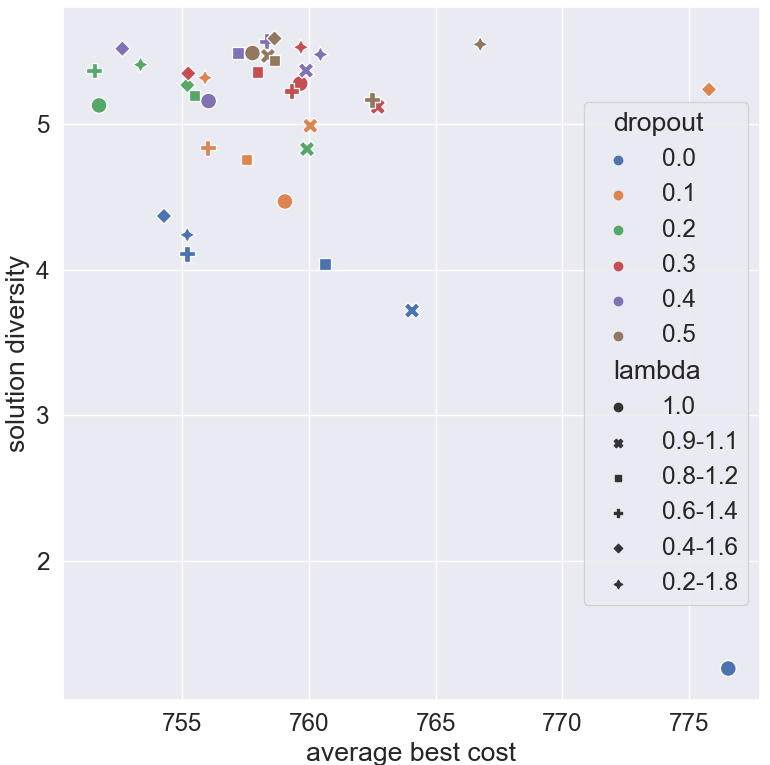}
\caption{Average cost and solution diversity of the solution after local search depending on the route shape parameter $\lambda$ and the dropout parameter $p_d$}
\label{FIG:hc}
\end{figure}

\begin{algorithm}
\caption{Local Search}
\label{alg2}
\begin{algorithmic}[1]
\STATE load inputs: vehicles, locations, chosen initial solution
\STATE data = vehicles, locations
\STATE select iterations limit $i_{ls}^{max}$

\FOR {$i_{ls}$ = 1 : $i_{ls}^{max}$}
    \STATE select strategy number $s$ based on $i_{ls}$
    \IF {$s$ == 1 or $s$ == 2}
        \STATE select number of combinations to test $g_c$ 
        \STATE select number of groups to test $g_n$
        \STATE select group size $g_s$
    \ENDIF
    \IF {$s$ == 1}
        \STATE solution = find\_best\_swap\_vertices(data, solution, $g_c$, $g_n$, $g_s$)
    \ELSIF{$s$ == 2}
        \STATE solution = find\_best\_move\_vertices(data, solution, $g_n$, $g_s$)
    \ELSIF{$s$ == 3}
        \STATE solution = find\_best\_swap\_drives(data, solution)
    \ELSIF{$s$ == 4}
        \STATE solution = find\_best\_move\_drive(data, solution)
    \ELSIF{$s$ == 5}
        \STATE solution = find\_best\_directions(data, solution)
    \ENDIF
\ENDFOR

\WHILE{solution is the same as in iteration before}
    \FOR {$s$ = 1 : 6}
        \STATE REPEAT 6-21 with max value of $g_c$ 
    \ENDFOR
\ENDWHILE

\RETURN local optimum solution

\end{algorithmic}
\end{algorithm}

\subsection{Adaptive memory}

The concept of adaptive memory is originally proposed in \cite{b18} to solve the VRP with time windows. Adaptive memory procedure (AMP) takes solutions obtained by multiple runs of a different algorithm as diversification of the search process. Here, we use local search and a certain percentage of the best solutions is selected. Percentage of "good" solutions among all obtained solutions is determined by $p_s$ as the retention parameter. Other solutions are ignored in further execution. Further on, common parts (a series of several vertices) of "good" solutions are identified, and transformed into a new vertex. This way, the problem is simplified by reducing the total number of vertices, and further search is performed on the reduced set. The demand of a new vertex is equal to the sum of the demands of the vertices from which it is created. The new time windows are the intersection of the time windows of the vertices it contains, and skills are the intersection of the skills of the vertices it contains. After the problem is reduced, the procedure is repeated for an arbitrary number of iterations. Finally, the solution with the lowest value of the cost function is selected. The pseudo-code is shown in Algorithm \ref{alg3}.

\begin{algorithm}
\caption{Adaptive Memory Procedure}
\label{alg3}
\begin{algorithmic}[1]
\STATE data $=$ vehicles, locations

\FOR {i\_1 $=$ 1 : iteration\_limit\_1}
\STATE select new iteration\_limit\_2
\STATE select good solutions percentage $p_s$
\STATE initialize solutions\_list
\FOR {i\_2 $=$ 1 : iteration\_limit\_2}

\STATE initial\_solution = run\_RCCW(data)
\STATE local\_optimum\_solution = run\_local\_search(data, initial\_solution)
\STATE append local\_optimum\_solution to solutions\_list
\ENDFOR
\IF {i\_1 $<$ iteration\_limit\_1}
\STATE keep best $p_s$ of solutions\_list
\STATE vertices $=$ select\_vertices\_to\_connect(solutions\_list)
\STATE data $=$ reduce\_data(data, vertices)
\ELSE
\STATE best\_solution = select\_best\_solutions(solution\_list, 1)
\ENDIF
\ENDFOR

\end{algorithmic}
\end{algorithm}

Also, the cost after each iteration of AMP is analyzed in Table~II. The analysis was carried out on 5 instances from \cite{b44}. Each iteration consists of 15 executions of RCCW and local search, and the average cost and standard deviation are shown per iteration. The standard deviation of the cost decreases through the iterations, mostly due to the locked "good" parts of the solution that do not change in further optimization. This leaves the possibility of searching only those parts of the graph that were poorly optimized until then. This also results in a lower average cost, which is the final goal of this algorithm. In the last iterations, the cost stagnates and the search stops after the fourth iteration. At the end, a more detailed 3-opt local search is performed on the best solution of the fourth iteration (by searching the reduced problem).

\begin{table}[htbp]
\caption{Cost analysis through AMP iterations}
\label{tab2}
\begin{center}
\begin{tabular}{|l|l|l|l|l|l|l|}
\hline
\scalebox{0.93}{Iteration} & value & CMT1 & CMT2 & CMT3 & \scalebox{0.93}{CMT11} & \scalebox{0.93}{CMT12} \\
\hline
1 & avg. & 583.62 & 902.89 & 923.36 & 1259.98 & 936.63 \\
& std. & 39.64 & 23.47 & 53.94 & 172.60 & 77.87 \\
\hline
2 & avg. & 566.63. & 890.20 & 898.07 & 1121.10 & 876.87 \\
& std. & 16.78 & 18.24 & 15.59 & 50.76 & 31.70 \\
\hline
3 & avg. & 561.55 & 881.30 & 886.72 & 1074.16 & 875.54 \\
& std. & 13.55 & 13.75 & 16.17 & 8.29 & 19.86 \\
\hline
4 & avg. & 557.69 & 868.86 & 873.35 & 1045.83 & 849.96 \\
& std. & 0.15 & 5.71 & 0.00 & 0.88 & 2.53 \\
\hline
5 & avg. & 533.48 & 857.99 & 855.19 & 1044.65 & 822.82 \\
\hline
\end{tabular}
\end{center}
\end{table}

\section{Benchmarks}

The AMP is tested on a last-mile delivery dataset provided by an industrial partner. Since the optimal solutions are not known for problem instances in this dataset, we added testing on one of the standardized benchmarks. We chose the benchmark from \cite{b44} to compare the AMP with the current state of art.

\subsection{Standardized benchmark}
Benchmark from \cite{b44} is one of the most commonly used in the field of VRP. It consists of 14 problem instances, 7 of which refer to CVRP, and the other 7 add time windows to vehicles and unloading time at each of the locations (CVRP+VRPTW). Delivery points are divided into city-like clusters in 4 instances, while in the other 10 instances their locations are completely or mostly randomized.

\subsection{Real world benchmark}
After extensive search we were unable to find a benchmark that includes all the considered VRP variants for real world application. Therefore, a real case study with data supplied from the industry partner is considered. The partner company delivers fresh products in vicinity of one of the major Croatian coastal cities. The studied data involves 118 distinct locations covering area of about 3350 $\text{km}^2$. The data consists of list of locations with their coordinates, time windows (working hours) defining when is possible to deliver orders, and packages included in clients order. Furthermore, each package has defined volume and mass that are needed for defining the remaining capacity of a vehicle in each step of delivery route. Distances and travel times are taken from navigation service and corresponding matrices are generated and used as part of this benchmark. Traveling by ferry to island locations is taken into account in the matrices as well. 

Calculating travel times for all pairs of locations requires over 1 million calculations ($I \cdot I \cdot K$). Therefore, for creation of variable travel times matrix, delivery points are empirically clusterized in 15 groups (clusters) with regard to their geographical locations: 8 groups in the city and suburbs and 7 in neighboring towns. Travel times between cluster centers are generated with sample of 15 minutes during a day. These travel times are normalized and used as time profiles (scalers) for travel times between locations. The same time profile is used for any pair of locations between two distinct clusters. Within a cluster, each pair has its own time profile. This approach assumes similar profile of time change for travel duration between locations of different clusters due to the same roads taken and traffic jam along them. Typical time profiles are shown in Fig.~6. The difference in driving time with static time matrix and rush-hour peaks with variable time matrix ranges from 3\% to 154\%.

\begin{figure}[!h]
\centering
\includegraphics[width=.48\textwidth]{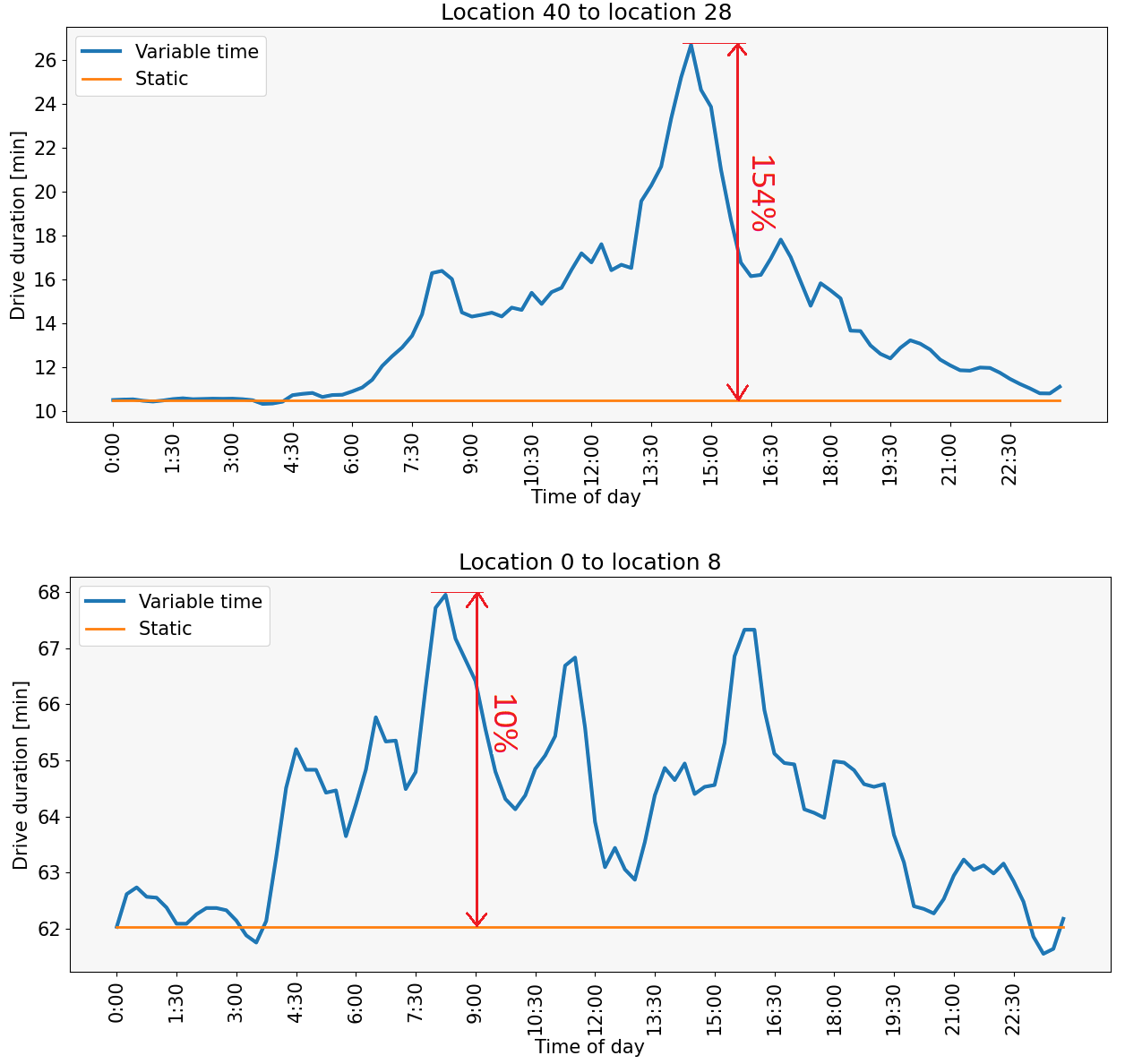}
\caption{Comparison of driving times between two locations with static and dynamic time matrices for a real-world scenario.}
\label{FIG:hc}
\end{figure}

The test set contains orders during five different days. The number of visited locations per day vary from 61 to 106, with average of 83 locations per day. Number of delivered packages per locations ranges from 1 to 31 packages, and with 330 to 684 packages per day. Finally, entirety of all orders consists of more than 2300 delivered packages over five days with 418 location visits.

\section{Results}

The proposed algorithm is tested on standardized and real world benchmarks. Since the AMP includes stochastic elements, we tested it on 10 runs and recorded the best and average case. For comparison, the currently best found solution for the chosen standardized benchmark, as reported in \cite{b51}, is used.  In addition, two algorithms from OR-Tools open source software for combinatorial optimization \cite{b49} are added. For the first OR-Tools solution, conventional Clarke-Wright with tabu search is chosen. It consists of searching for an initial solution using the Clarke-Wright algorithm, after which tabu search is used as a metaheuristic. Tabu search is selected as one of the most used methods due to inability to select AMP. In the sequel, we use abbreviation TS for the combination of Clarke-Wright algorithm and tabu search. For the second OR-Tools approach, the path cheapest arc method is chosen for the initial solution and Guided Local Search (GLS) as a metaheuristic. Abbreviations and descriptions of algorithms are shown in Table~III.
The OR-Tools is currently a widely used software. It can solve problems with a much larger number of constraints than the algorithms used to obtain the optimal solutions from the standardized benchmark. One of the variants that this software does not solve, and which together with our industrial partner we consider one of the most important, is TDVRP.

\begin{table}[htbp]
\caption{List of tested algorithms with short descriptions}
\label{tab2}
\begin{center}
\begin{tabular}{|l|l|l|}
\hline
Abbreviation & Applied algorithms & source \\
\hline
& Randomized Clarke-Wright + & proposed in \\
AMP & Local search + & this paper \\
& Adaptive memory & \\
\hline
TS & Clarke-Wright + & open-source \\
& Tabu search & from OR-Tools \\
\hline
GLS & Path cheapest arc + & open-source \\
& Guided local search & from OR-Tools \\
\hline
\end{tabular}
\end{center}
\end{table}

\begin{table*}[htbp]
\caption{Comparison of different approaches on the standardized benchmark from \cite{b44}}
\begin{center}
\label{tab3}
\begin{tabular}{|l|c|c|r|r|r|r|r|r|c|}
\hline
Instance & Size & Variants & \multicolumn{1}{|c|}{Best known} & GLS cost & TS cost & \multicolumn{2}{|c|}{AMP cost} & time [s]& AMP \\
\cline{7-8}
&&& \multicolumn{1}{|c|}{cost} & & & best run & average run & & suboptimality \\
\hline
CMT1 &50&CVRP & 524.61 \cite{b17} & 531.02 & 524.61 & 524.93 & 532.54 & 159 & 0.06\% \\
CMT2 &75&CVRP& 835.26 \cite{b17} & 849.19 & 854.17 & 846.83 & 853.80 & 336 & 1.39\% \\
CMT3 &100&CVRP& 826.14 \cite{b17} & 836.18 & 834.56 & 838.81 & 847.47 & 648 & 1.53\% \\
CMT4 &150&CVRP& 1,028.42 \cite{b17} & 1,043.62 & 1,075.38 & 1,065.81 & 1,075.16 & 1,743 & 3.64\% \\
CMT5 &199&CVRP& 1,291.29 \cite{b50} & 1,330.52 & 1,345.49 & 1,349.72 & 1,358.58 & 4,050 & 4.52\% \\
CMT6 &50&CVRP+VRPTW& 555.43 \cite{b17} & 558.99 & 556.68 & 555.43 & 560.81 & 142 & 0.00\% \\
CMT7 &75&CVRP+VRPTW& 909.68 \cite{b17} & 918.27 & 932.15 & 925.43 & 936.24 & 345 & 1.73\% \\
CMT8 &100&CVRP+VRPTW& 865.95 \cite{b17} & 867.14 & 870.32 & 887.34 & 892.57 & 598 & 2.47\% \\
CMT9 &150&CVRP+VRPTW& 1,162.55 \cite{b17} & 1,183.26 & 1,189.36 & 1,204.92 & 1,214.22 & 1,942 & 3.38\% \\
CMT10 &199&CVRP+VRPTW& 1,395.85 \cite{b18} & 1,440.07 & 1,440.24 & 1,456.92 & 1,469.21 & 4,291 & 4.38\% \\
CMT11 &120&CVRP& 1,042.12 \cite{b17} & 1,051.49 & 1,177.10 & 1,043.89 & 1,046.57 & 1,061 & 0.17\% \\
CMT12 &100&CVRP& 819.56 \cite{b17} & 819.56 & 825.95 & 821.69 & 829.01 & 586 & 0.26\% \\
CMT13 &120&CVRP+VRPTW& 1,541.14 \cite{b17} & 1,580.78 & 1,614.62 & 1,556.83 & 1,560.55 & 1,087 & 1.02\% \\
CMT14 &100&CVRP+VRPTW& 866.37 \cite{b17} & 908.52 & 924.15 & 868.42 & 875.26 & 542 & 0.24\% \\
\hline
average &&& 976.03 & 994.19 & 1,011.75 & 996.21 & 1,003.71 && 1.77\% \\
\hline

\end{tabular}
\end{center}
\end{table*}

\subsection{Results on the standardized benchmark}

In the standardized benchmark test, we compared the AMP with OR-Tools algorithms and best known solutions on two variants: CVRP and VRPTW. The test results are given in Table~IV. The average cost for the best known solutions on 14 instances is 976.03. The AMP, with the average of the best solutions of 996.21, proved to be better than TS algorithm. This is partly due to the dropout, which increased the search space for solutions, and partly due to different metaheuristics used. The average cost of AMP for the best of 10 runs is 996.21, which is slightly worse than GLS algorithm. Column time refers to the execution time of the AMP, and the same time is given to TS and GLS algorithms. A solution equal to the best known was found for instance CMT6 (from \cite{b44}), and the average suboptimality on the set of 14 instances is 1.77\%.
The value for suboptimality of each instance is calculated from the ratio of AMP cost and best known cost. It should be noted that the suboptimality is lower at instances with fewer delivery points, which is expected because, in that case, a bigger part of the hyperplane is searched.

\subsection{Results on the real world benchmark}

We performed the second test on industrial partner data containing multiple variants including TDVRP. We compared the proposed AMP that solves all the variants listed in the paper with GLS algorithm. Testing was performed on listed all variants except TDVRP which the tested GLS implementation cannot solve.

The obtained results are shown in Table~V. In the table, the columns denote, respectively: day number, number of packages, stores and vehicles, and finally, cost in seconds for GLS and AMP without inclusion of TDVRP. For the case with time matrices, the last time interval (6:45 a.m.) before loading packages into the vehicles (7:00 a.m.) is selected.
Finally, instead of static time matrix, dynamic ones (TDVRP) were added. The GLS solution remained the same, but a new cost is calculated that increased due to traffic jams (increased $c_{nijk}$) within the interval [8:00-14:30] from (3) and (4). Due to traffic congestions, the soft time window $ e^{vrs}_i, \forall i$ from equation (5) is exceeded. Therefore the 'GLS exceeded' column is added, which shows how long delay was caused. The AMP found a new solution without exceeded time. It is shown in the AMP cost column.

\begin{table*}[htbp]
\caption{Comparison of different algorithms in real world benchmark, delivery time cost}
\label{tab4}
\begin{center}
\begin{tabular}{|c|c|c|c|r|r|r|r|r|}
\hline

Day & Num. of & Num. of & Num. of & \multicolumn{2}{|c|}{static time matrix} & \multicolumn{3}{|c|}{variable time matrix (TDVRP)} \\

\cline{5-9}
no. & packages & stores & vehicles & GLS cost & AMP cost & GLS cost & GLS exceeded & AMP cost \\

\hline
1 & 335 & 61 & 4 & 42,145s & 42,149s & 43,771s & 0s & 42,263s \\
2 & 445 & 86 & 5 & 74,289s & 74,413s & 76,978s & 694s & 75,845s \\
3 & 437 & 82 & 5 & 61,142s & 61,296s & 63,685s & 1,959s & 62,758s \\
4 & 454 & 83 & 4 & 63,055s & 63,221s & 66,813s & 388s & 64,506s \\
5 & 704 & 106 & 5 & 69,316s & 70,075s & 71,606s & 0s & 70,914s \\
\hline
average & & & & 61,989s & 62,231s & 64,571s & 608s & 63,257s \\
\hline
\end{tabular}
\end{center}
\end{table*}

The average cost as first element of (11) without time windows for GLS is lower by 0.39\% with static matrix when comparing GLS and AMP. This cost refers to the total distance traveled for CVRP and to the total delivery time for CVRP+VRPTW. However, this changes with the inclusion of a variable time matrix. In this case, the driving time is 2.03\% shorter for AMP than for GLS solution. Also, the average delay for GLS is 608 seconds, which is 0.94\% of the total time, while with AMP all orders are delivered within working hours. In the dataset, most of the delivery points are outside the city center and the city itself has less than 200 thousand inhabitants. In the case of a larger city, the savings on GLS cost, and especially GLS exceeded, would be significantly higher.

\begin{table*}[htbp]
\caption{Comparison of different algorithms in real world benchmark, monetary cost}
\label{tab5}
\begin{center}
\begin{tabular}{|c|c|c|c|r|r|r|r|}
\hline
Day & Num. of & Num. of & Num. of & \multicolumn{2}{|c|}{static time matrix} & \multicolumn{2}{|c|}{variable time matrix (TDVRP)} \\
\cline{5-8}
no. & packages & stores & vehicles & GLS cost & AMP cost & GLS cost & AMP cost \\

\hline
1 & 335 & 61 & 4 & 51.229€ & 51.213€ & 52.193€ & 51.990€ \\
2 & 445 & 86 & 5 & 103.442€ & 103.344€ & 164.501€ & 104.824€ \\
3 & 437 & 82 & 5 & 70.596€ & 70.529€ & 100.003€ & 71.856€ \\
4 & 454 & 83 & 4 & 74.725€ & 74.662€ & 90.936€ & 75.806€ \\
5 & 704 & 106 & 5 & 84.503€ & 85.002€ & 85.859€ & 85.484€ \\
\hline
average & & & & 76.889€ & 76.950€ & 98.698€ & 77.992€ \\
\hline
\end{tabular}
\end{center}
\end{table*}

A more detailed case study with a monetary cost comparison is given in Table~VI, but with a monetary cost. The difference is that in this case GLS exceeded time is not observed, but the cost of exceeding the time window is immediately calculated as an additional cost in the GLS cost column. Although many parameters can be added to the monetary cost, only those dependent on the VRP solution for this dataset are analyzed here: driving cost characteristic of each arc and cost for deliveries outside the ideal time window. Therefore, the costs included in the GLS cost column are covered by (10). Base cost $c_{nijk}$ is calculated as the product of the fuel consumed (2) and the fuel price. An additional cost is calculated for such a delay according to (8). The parameters used are $a^e_0=1.5$, $a^e_1=0.001$ and $a^e_2=0.0007$. 
This means that each delay brings a cost of 1.5€, each second of delay costs 0.001€, and each second of delayed unloading costs 0.0007€. A similar additional cost can be calculated for the delay of the vehicle when returning to the depot according to equation (6) with the interval [8:00-13:30], but in this case no such situation occurred. At the delivery points in this example, late deliveries are possible, but their cost is relatively high compared to the cost of fuel consumption. Because of this, on days 2, 3 and 4, the GLS cost with variable time matrix increased significantly.

Using static time matrix to optimize the monetary cost, GLS found a solution with 0.08\% lower cost than AMP, similarly as for delivery time cost in Table~V. However, significant savings were achieved with the variable time matrix application, mostly due to late deliveries in the case of GLS. Solutions for days 1 and 5 have no late deliveries and the costs are similar, but in days 2, 3 and 4, the cost of GLS is significantly higher because of the late deliveries. Averaged over all five days, AMP achieved savings of 20.98\% compared to GLS.

The code was written in the Python programming language, and testing is performed on one core of Intel i5-10210U CPU (1.6GHz). The execution time is equal for GLS and AMP in Tables IV and V, and is between 233 and 814 seconds.

\section{Conclusion}

This paper describes an algorithm for solving a combination of different variants of VRP: CVRP, VRPTW, VRPSTW, HVRP, MTVRP, SVRP, SDVRP, TDVRP. 
The Adaptive Memory Procedure (AMP) metaheuristic is used, which searches for best common parts of solutions obtained by multiple runs of simpler methods, thereby reducing the space for further searches. To find the initial solution, Randomized Constrained Clarke-Wright algorithm is proposed. Such obtained solution is improved by classical methods of local search. Standardized benchmark testing show an average suboptimality of 1.77\% on CVRP and VRPTW problems. A combination of all the above methods is tested on real world dataset. The AMP achieved 2.03\% savings in delivery time and 20.98\% on the variable part of the monetary cost compared to the current state of the art. The results show promising contribution of high practical and economic significance when considering real-world problems, and a potential to further reduce the suboptimality.

\vspace{12pt}

\end{document}